\documentclass[11pt]{amsart}
\usepackage{latexsym}
\usepackage{amsmath}
\usepackage{amssymb}

\newcommand{\eproof}{\mbox{\ }\hfill $\Box$ \par \vskip 10pt}

\newtheorem{Theorem}{Theorem}[section]
\newtheorem{lemma}[Theorem]{Lemma}
\newtheorem{prop}[Theorem]{Proposition}

\numberwithin{equation}{section}

\def\cal{\mathcal}

\begin{document}

\title[Approximation of the elastic Dirichlet-to-Neumann map]{Approximation of the elastic Dirichlet-to-Neumann map}

\author[G. Vodev]{Georgi Vodev}

\address {Universit\'e de Nantes, Laboratoire de Math\'ematiques Jean Leray, 2 rue de la Houssini\`ere, BP 92208, 44322 Nantes Cedex 03, France}
\email{Georgi.Vodev@univ-nantes.fr}

\date{}

\begin{abstract} We study the Dirichlet-to-Neumann map for the stationary linear equation of elasticity in a bounded domain
in $\mathbb{R}^d$, $d\ge 2$, with smooth boundary.
We show that it can be approximated by a pseudodifferential operator on the boundary with a matrix-valued symbol
and we compute the principal symbol modulo conjugation by unitary matrices. 
\quad

\end{abstract} 

\maketitle

\setcounter{section}{0}
\section{Introduction}

Let $\Omega\subset\mathbb{R}^d$, $d\ge 2$, be a bounded, connected domain with a $C^\infty$ smooth boundary $\Gamma=\partial\Omega$, and consider
the stationary isotropic linear equation of elasticity

\begin{equation}\label{eq:1.1}
\left\{
\begin{array}{l}
(\Delta_{\lambda,\mu}+\tau^2n(x))u=0\quad \mbox{in}\quad\Omega,\\
u=f\quad\mbox{on}\quad\Gamma,
\end{array}
\right.
\end{equation}
where $\tau\in \mathbb{C}$, ${\rm Re}\,\tau>0$, $|\tau|\gg 1$, $u=(u_1,...,u_d)$, $f=(f_1,...,f_d)$, 
and $\Delta_{\lambda,\mu}$ denotes the elastic Laplacian defined by 
$$\left(\Delta_{\lambda,\mu}u\right)_i=\sum_{j=1}^d\partial_{x_j}\left(\sigma_{ij}(u)\right),\quad i=1,...,d,$$
where $$\sigma_{ij}(u)=\lambda{\rm div}\, u\delta_{ij}+\mu\left(\frac{\partial u_i}{\partial x_j}+\frac{\partial u_j}{\partial x_i}\right)$$
is the stress tensor, $\delta_{ij}=1$ if $i=j$, $\delta_{ij}=0$ if $i\neq j$. 
Here $\lambda,\mu\in C^\infty(\overline{\Omega})$ are scalar real-valued functions called Lam\'e parameters supposed to satisfy the condition
\begin{equation}\label{eq:1.2}
\mu(x)>0,\quad\lambda(x)+\mu(x)>0,\quad\forall x\in \overline{\Omega}.
\end{equation}
The scalar function $n\in C^\infty(\overline{\Omega})$ in (\ref{eq:1.1}) is the density and is supposed to be strictly positive. 
It is easy to see that the elastic Laplacian can be written in the form
$$\Delta_{\lambda,\mu}u=\mu\Delta u+(\lambda+\mu)\nabla(\nabla\cdot u)$$
modulo a first-order matrix-valued differential operator, where $\Delta$ and $\nabla$ denote the Euclidean Laplacian and gradient, respectively. 

The natural Neumann boundary condition for the elastic equation is $B_{\lambda,\mu}u=0$, where
$$\left(B_{\lambda,\mu}u\right)_i=\sum_{j=1}^d\sigma_{ij}(u)\nu_j,\quad i=1,...,d,$$
$\nu=(\nu_1,...,\nu_d)$ being the Euclidean unit normal to $\Gamma$. 
 We define the elastic Dirichlet-to-Neumann map
$$N(\tau):H^1(\Gamma;\mathbb{C}^d)\to L^2(\Gamma;\mathbb{C}^d)$$
by
$$N(\tau)f=B_{\lambda,\mu}u|_{\Gamma}$$
where $u$ and $f$ satisfy the equation (\ref{eq:1.1}).

The equation (\ref{eq:1.1}) describes the propagation of elastic waves in
$\Omega$ with a frequency $\tau$. It is well-known that the elastic waves are superpositions of two waves,
called S and P waves, mooving with speeds $\sqrt{\frac{\mu}{n}}$ and $\sqrt{\frac{2\mu+\lambda}{n}}$, respectively. 
From purely mathematical point of view, this is explained by the fact that the principal symbol, $P$, of the operator
$-\Delta_{\lambda,\mu}$ can be decomposed as 
$$P(x,\xi)=c_s(x)\Pi_s(\xi)+c_p(x)\Pi_p(\xi)$$
where $c_s=\mu$, $c_p=2\mu+\lambda$, $\Pi_s(\xi)+\Pi_p(\xi)=\xi^2I_d$, $I_d$ being the identity $d\times d$ matrix, and 
$\Pi_p(\xi)=\xi\otimes\xi$. Throughout this paper, given two vectors $\xi=(\xi_1,...,\xi_d)\in \mathbb{C}^d$, 
$\eta=(\eta_1,...,\eta_d)\in \mathbb{C}^d$, we will denote by $\xi\otimes\eta$ 
the matrix defined by
$$(\xi\otimes\eta)g=\langle \xi,g\rangle\eta,\quad g\in \mathbb{C}^d.$$
Hereafter, $\langle \xi,g\rangle:=\xi_1g_1+...+\xi_dg_d$ and $\xi^2:=\langle \xi,\xi\rangle$. 

Note that the existance of two different speeds implies that the boundary value problem (\ref{eq:1.1}) has two disjoint
glancing regions $\Sigma_s$ and $\Sigma_p$ defined by
$$\Sigma_s=\left\{(x',\xi')\in T^*\Gamma:c_{s,0}(x')r_0(x',\xi')-n_0(x')=0\right\},$$
$$\Sigma_p=\left\{(x',\xi')\in T^*\Gamma:c_{p,0}(x')r_0(x',\xi')-n_0(x')=0\right\},$$
where $c_{s,0}=c_s|_{\Gamma}$, $c_{p,0}=c_p|_{\Gamma}$, $n_0=n|_{\Gamma}$, and $r_0\ge 0$ is 
the principal symbol of the operator $-\Delta_{\Gamma}$. Here
$\Delta_{\Gamma}$ denotes the negative Laplace-Beltrami operator on $\Gamma$ with Riemannian metric induced by the Euclidean one. 
Set $h=({\rm Re}\,\tau)^{-1}$ if 
${\rm Re}\,\tau\ge |{\rm Im}\,\tau|$ and $h=|{\rm Im}\,\tau|^{-1}$ if 
$|{\rm Im}\,\tau|\ge {\rm Re}\,\tau$, $z=h\tau$ and $\theta=|{\rm Im}\,z|\le 1$. 
Clearly, in the first case we have $z=1+i\theta$, while in the second case we have $\theta=1$.  When $\theta>0$ we introduce the functions 
$$\rho_s(x',\xi',z)=\sqrt{-r_0(x',\xi')+z^2n_0(x')/c_{s,0}(x')},\quad {\rm Im}\,\rho_s>0,$$
$$\rho_p(x',\xi',z)=\sqrt{-r_0(x',\xi')+z^2n_0(x')/c_{p,0}(x')},\quad {\rm Im}\,\rho_p>0.$$
Our goal in the present paper is to approximate
the operator 
$${\cal N}(z,h):=-ihN(z/h)$$
 by a matrix-valued $h-\Psi$DO similarly to the Dirichlet-to-Neumann operator associated to
the Helmholtz equation (see \cite{kn:V1}, \cite{kn:V2}) or that one associated to the Maxwell equation (see \cite{kn:V3}).  
We also compute the principal symbol in terms of the functions $\rho_s$, $\rho_p$ (see Lemma 5.5). Denote by $M_2$ the $2\times 2$
matrix with entries $M_{ij}$ given by
$$M_{11}=\frac{z^2n_0\rho_s}{r_0+\rho_s\rho_p},\quad M_{22}=\frac{z^2n_0\rho_p}{r_0+\rho_s\rho_p},$$
$$-M_{21}=M_{12}=-2\mu_0\sqrt{r_0}+\frac{z^2n_0\sqrt{r_0}}{r_0+\rho_s\rho_p},$$
where $\mu_0=\mu|_{\Gamma}=c_{s,0}$. 
When $d\ge 3$ we set
$$M_d=\widetilde M_2+\mu_0\rho_s(I_d-\widetilde{I_2}).$$
Throughout this paper, given a $2\times 2$
matrix $M$ with entries $M_{ij}$, we denote by $\widetilde M$ the $d\times d$ matrix with entries 
$\widetilde M_{ij}=M_{ij}$ if $1\le i,j\le 2$, $\widetilde M_{ij}=0$ otherwise.  Given a partition of the unity 
$\kappa_\ell\in C^\infty(T^*\Gamma\setminus 0)$, $0\le\kappa_\ell\le 1$, $\ell=1,...,L$, $\sum_{\ell=1}^L\kappa_\ell=1$, introduce the function
\begin{equation}\label{eq:1.3}
m_d=\sum_{\ell=1}^L\kappa_\ell J_\ell M_dJ_\ell^{-1}
\end{equation}
where $J_\ell(x',\xi')\in C^\infty(T^*\Gamma\setminus 0)$ are matrix-valued functions, homogeneous of order zero in $\xi'$, and such that $J_\ell^{-1}=J_\ell^t$.
Our main result is the following

\begin{Theorem} 
Let $\theta\ge h^{2/5-\epsilon}$ and $0<h\ll 1$, where $0<\epsilon\ll 1$ is arbitrary. Then for every $f\in H^3(\Gamma;\mathbb{C}^d)$ 
we have the estimate
\begin{equation}\label{eq:1.4}
\left\|{\cal N}(z,h)f-{\rm Op}_h(m_d)f\right\|_{L^2(\Gamma;\mathbb{C}^d)}\lesssim h\theta^{-2}\left(1+(d-2)\theta^{-1/2}\right)\|f\|_{H_h^3(\Gamma;\mathbb{C}^d)}
\end{equation}
 where $m_d\in C^\infty(T^*\Gamma)$ is of the form (\ref{eq:1.3}) with a suitable partition of the unity $\kappa_\ell$ and 
  matrix-valued functions $J_\ell$ independent of $\lambda$, $\mu$ and $n$. 
  When $d=2$ the functions $J_\ell$ do not depend on the variable $\xi'$. 
 \end{Theorem}

Hereafter the Sobolev spaces are equipped with the $h$-semiclassical norm.
Note that much better estimates for the Dirichlet-to-Neumann operator associated to
the Helmholtz equation are proved in \cite{kn:V1}, \cite{kn:V2}. This is due to the fact that one can construct a much better
parametrix near the boundary for the Helmholtz equation than that one for the equation (\ref{eq:1.1}) we construct in the present paper.
Indeed, such a parametrix is built in \cite{kn:V1}, \cite{kn:V2} in the form of an oscilatory integral with a complex-valued 
phase function and an 
amplitude satisfying the eikonal and transport equations mod ${\cal O}(x_1^N)$, respectively,  where $N\gg 1$
is arbitrary and $0<x_1\ll 1$ denotes the normal variable near the boundary,
that is, the distance to $\Gamma$.
 Thus the parametrix satisfies the Helmholtz equation
modulo an error term which is given by an oscilatory integral with amplitude of the form ${\cal O}(x_1^N)+{\cal O}(h^N)$. 
In the case of the equation (\ref{eq:1.1}), however, it is very hard to solve the transport equations, especially when the boundary
data $f$ is microlocally supported in a neighbourhood of the glancing regions. That is why we build in the present paper a less
accurate parametrix for the equation (\ref{eq:1.1}) which does not require to solve transport equations. 
In this case the parametrix is a sum of two oscilatory integrals with two complex-valued 
phase functions corresponding to the two speeds of propagation of the elastic waves. Each of these  phase functions
satisfies the same eikonal equation as in the case of the Helmholtz equation solved in \cite{kn:V1}. 
The parametrix satisfies the equation (\ref{eq:1.1}) 
modulo an error term which is given by a sum of two oscilatory integrals with amplitudes of the form ${\cal O}(x_1^N)+{\cal O}(h)$. 
To estimate the
difference between the exact solution to equation (\ref{eq:1.1}) and its parametrix we use the a priori estimate (\ref{eq:4.2}).
Most probably, the estimate (\ref{eq:1.4}) is not optimal and could be improved if one manages to build a better parametrix. In particular,
for some applications it is better to have the $L^2$ norm in the right-hand side of (\ref{eq:1.4}) instead of the $H_h^3$ one. 
To do so, one needs to construct a better parametrix in the deep elliptic region, only,
that is, in $\{r_0\gg 1\}$. Recall that the approximation of the Dirichlet-to-Neumann map is usually used to get parabolic regions free of transmission eigenvalues (see \cite{kn:V1}, \cite{kn:V2}, \cite{kn:V3}). 

Note that microlocal parametrices have been recently constructed in \cite{kn:HINZ}, \cite{kn:SUV}, \cite{kn:Z} for the wave elastic
equation in the case $d=3$. All these parametrices, however, are very different from the parametrix we construct in the present paper.
In particular, they are not valid near the glancing regions. In contrast, our parametrix remains valid even when the boundary data
is microlocally supported near $\Sigma_s$ and $\Sigma_p$, provided $\theta\ge h^{2/5-\epsilon}$. In \cite{kn:Z} the principal symbol
of the DN map has been computed explicitly still in the context of the wave elastic equation and $d=3$. Note that the formula in
\cite{kn:Z} agrees with that one we get in the present paper modulo a conjugation by 
a unitary matrix and after making a suitable change in the notations. 
In \cite{kn:SV1} a full parametrix was constructed for the stationary elastic equation in the exterior of a strictly convex body when 
$d=3$ and the Lam\'e parameters being constants. In this case the parametrix for the elastic equation can be expressed in terms of the
parametrix for the Helmholtz equation, which in turn is well-known. Similarly, one can construct a parametrix in the elliptic region 
for the stationary elastic equation in the exterior of an arbitrary compact body (see \cite{kn:SV2}).
 
 \section{Preliminaries}
 
 Throuhout this paper we will denote by $e_1,e_2,...,e_d\in \mathbb{R}^d$ the vectors $(1,0,...,0)$,
$(0,1,...,0)$, ..., $(0,0,...,1)$, respectively. Given $\xi=(\xi_1,\xi_2,...,\xi_d)\in \mathbb{C}^d$, introduce the $d\times d$ matrix
$$U_0(\xi)=\xi_1I_d+\sum_{j=2}^d\xi_j(e_j\otimes e_1-e_1\otimes e_j).$$
Set $\xi^2=\sum_{j=1}^d\xi_j^2$ and $|\xi|^2=\sum_{j=1}^d|\xi_j|^2$.  
 In this section we will prove the next two lemmas.
 
\begin{lemma} The matrix $U_0$ satisfies
\begin{equation}\label{eq:2.1}
U_0(\xi)\xi=\xi^2e_1.
\end{equation}
Moreover, the matrices $e_1\otimes e_1$ and $Z_0(\xi):=U_0(\xi)U^t_0(\xi)$ commute. 
\end{lemma}
 
 {\it Proof.} We have
 $$U_0(\xi)\xi=\xi_1\xi+\sum_{j=2}^d\xi_j(\langle e_j,\xi\rangle e_1-\langle e_1,\xi\rangle e_j)$$
 $$=\sum_{j=1}^d\xi_1\xi_je_j+\sum_{j=2}^d\xi_j(\xi_j e_1-\xi_1e_j)=\sum_{j=1}^d\xi_j^2e_1$$
 which proves (\ref{eq:2.1}). 
 Set $\widetilde\xi=(0,\xi_2,...,\xi_d)$. Then $U^t_0(\widetilde\xi)=-U_0(\widetilde\xi)$ and hence 
 $$Z_0(\xi)=\xi_1^2I_d-U_0(\widetilde\xi)^2.$$ 
 On the other hand, it is easy to see that
 $$U_0(\widetilde\xi)(e_1\otimes e_1)=-\sum_{j=2}^d\xi_je_1\otimes e_j,$$
  $$(e_1\otimes e_1)U_0(\widetilde\xi)=\sum_{j=2}^d\xi_je_j\otimes e_1.$$
  Thus we get
  $$(e_1\otimes e_1)U_0(\widetilde\xi)^2=-\widetilde\xi^2e_1\otimes e_1=U_0(\widetilde\xi)^2(e_1\otimes e_1),$$
  and hence $Z_0(\xi)$ and $e_1\otimes e_1$ commute. 
  \eproof
  
  Given a matrix $A$ with entries $A_{ij}$, define its norm by $\|A\|=\max_{ij}\left|A_{ij}\right|$.
  
  \begin{lemma} Given any $\eta_1\in \mathbb{C}$ we have the formula
\begin{equation}\label{eq:2.2}
{\rm det}(U_0(\xi)+\eta_1e_1\otimes e_1)=(\xi^2+\xi_1\eta_1)\xi_1^{d-2}.
\end{equation}
 When $d\ge 3$ suppose that $\xi_1\in \mathbb{C}$, $\xi_1\neq 0$, and $\xi_k\in \mathbb{R}$, $k=2,...,d$. Then, if $\xi^2+\xi_1\eta_1\neq 0$,
 we have the estimate 
\begin{equation}\label{eq:2.3}
\left\|(U_0(\xi)+\eta_1e_1\otimes e_1)^{-1}\right\|\lesssim (|\xi|+|\eta_1|)|\xi^2+\xi_1\eta_1|^{-1}+(d-2)|\xi_1|^{-1}.
\end{equation}
\end{lemma}
  
  {\it Proof.} Denote the matrix $U_0(\xi)+\eta_1e_1\otimes e_1$ by ${\cal M}_d(\xi_1,...,\xi_d,\eta_1)$. The lemma is very easy to prove 
  when $d=2$. Indeed, in this case we have
  $${\cal M}_2(\xi_1,\xi_2,\eta_1)=
  \left(
\begin{array}{cc}
\xi_1+\eta_1&\xi_2\\
-\xi_2&\xi_1
\end{array}
\right)$$
  and hence ${\rm det}{\cal M}_2=\xi^2+\xi_1\eta_1$. Moreover, if $\xi^2+\xi_1\eta_1\neq 0$, we have
  $${\cal M}_2^{-1}(\xi_1,\xi_2,\eta_1)=(\xi^2+\xi_1\eta_1)^{-1}
  \left(
\begin{array}{cc}
\xi_1&-\xi_2\\
\xi_2&\xi_1+\eta_1
\end{array}
\right)$$
  and (\ref{eq:2.3}) in this case is obvious. When $d\ge 3$ the formula (\ref{eq:2.2}) can be proved by induction. Indeed, we have
  $${\rm det}{\cal M}_d(\xi_1,...,\xi_d,\eta_1)=\xi_1{\rm det}{\cal M}_{d-1}(\xi_1,...,\xi_{d-1},\eta_1)+
  (-1)^d\xi_d{\rm det}{\cal P}_{d-1}(\xi_1,...,\xi_d),$$
  where ${\cal P}_{d-1}$ denotes the $(d-1)\times(d-1)$ matrix with lines $(\xi_2,...,\xi_d)$, $(0,\xi_1,0,...,0)$,
  ..., $(0,...,\xi_1,0)$. Hence 
  $${\rm det}{\cal P}_{d-1}(\xi_1,\xi_2,...,\xi_d)=\xi_1^{d-2}{\rm det}{\cal P}_{d-1}(1,\xi_2,...,\xi_d)$$
  $$=\xi_1^{d-2}{\rm det}{\cal P}_{d-1}(1,0,...,0,\xi_d)=(-1)^d\xi_d\xi_1^{d-2}.$$
  Thus we get
   $${\rm det}{\cal M}_d(\xi_1,...,\xi_d,\eta_1)=\xi_1{\rm det}{\cal M}_{d-1}(\xi_1,...,\xi_{d-1},\eta_1)+\xi_d^2\xi_1^{d-2}.$$
   It is clear now that if (\ref{eq:2.2}) holds for ${\rm det}{\cal M}_{d-1}$, it holds for ${\rm det}{\cal M}_d$, as well. Thus we conclude that
   (\ref{eq:2.2}) holds for all $d$. 
   
   The proof of (\ref{eq:2.3}) when $d\ge 3$ is more delicate. Set $\widetilde\xi=(0,\xi')\in\mathbb{R}^d$,
   where $\xi'=(\xi_2,...,\xi_d)\in\mathbb{R}^{d-1}$, and suppose that there exists a $d\times d$ matrix $\Theta(\xi')$,
   homogeneous of order zero, such that 
$\Theta(\xi')e_1=e_1$, $\Theta(\xi')\widetilde\xi=|\xi'|e_2$, $\Theta(\xi')^{-1}=\Theta^t(\xi')$ and $\|\Theta(\xi')\|\le C$ with a constant $C>0$ independent of $\xi'$. 
Then, given any $g\in\mathbb{C}^d$, we have
  $$U_0(\widetilde\xi)\Theta^t(\xi')g=\langle \Theta(\xi')\widetilde\xi,g\rangle e_1-\langle \Theta(\xi')e_1,g\rangle \widetilde\xi
  =|\xi'|\langle e_2,g\rangle e_1-\langle e_1,g\rangle \widetilde\xi.$$
  Hence
  $$\Theta(\xi')U_0(\widetilde\xi)\Theta^t(\xi')g=|\xi'|\langle e_2,g\rangle \Theta(\xi')e_1-\langle e_1,g\rangle \Theta(\xi')\widetilde\xi$$
  $$=|\xi'|\langle e_2,g\rangle e_1-|\xi'|\langle e_1,g\rangle e_2=U_0(|\xi'|e_2)g,$$
  which implies 
  $$
  U_0(\xi)=\Theta^t(\xi')U_0(\xi_1e_1+|\xi'|e_2)\Theta(\xi').$$
  Since $\Theta^t(\xi')(e_1\otimes e_1)\Theta(\xi')=(e_1\otimes e_1)$, we obtain
  \begin{equation}\label{eq:2.4}
  \left(U_0(\xi)+\eta_1e_1\otimes e_1\right)^{-1}=\Theta(\xi')\left(U_0(\xi_1e_1+|\xi'|e_2)+\eta_1e_1\otimes e_1\right)^{-1}\Theta^t(\xi').
  \end{equation}
  This implies 
  \begin{equation}\label{eq:2.5}
  \left\|\left(U_0(\xi)+\eta_1e_1\otimes e_1\right)^{-1}\right\|\lesssim\left\|\left(U_0(\xi_1e_1+|\xi'|e_2)+\eta_1e_1\otimes e_1\right)^{-1}\right\|.
  \end{equation}
  On the other hand,
  \begin{equation}\label{eq:2.6}
  \left(U_0(\xi_1e_1+|\xi'|e_2)+\eta_1e_1\otimes e_1\right)^{-1}=\widetilde{{\cal M}_2^{-1}}(\xi_1,|\xi'|,\eta_1)
  +\xi_1^{-1}\sum_{j=3}^de_j\otimes e_j.
  \end{equation}
  It follows from (\ref{eq:2.6}) that
   \begin{equation}\label{eq:2.7}
  \left\|\left(U_0(\xi_1e_1+|\xi'|e_2)+\eta_1e_1\otimes e_1\right)^{-1}\right\|\le 
  \left\|{\cal M}_2^{-1}(\xi_1,|\xi'|,\eta_1)\right\|+|\xi_1|^{-1}.
  \end{equation}
  Clearly, in this case (\ref{eq:2.3}) follows from (\ref{eq:2.5}) and (\ref{eq:2.7}).
  
  It remains to see that such a matrix $\Theta(\xi')$ exists. When $d=3$ it is easy to see that the matrix
  $$\Theta(\xi')=\left(
\begin{array}{ccc}
1&0&0\\
0&\xi_2/|\xi'|&\xi_3/|\xi'|\\
0&-\xi_3/|\xi'|&\xi_2/|\xi'|
\end{array}
\right)$$
  has the desired properties. When $d\ge 4$, however, it is hard (and probably impossible) to find only one global matrix $\Theta(\xi')$
  with these properties. Such a matrix, however, exists locally. Indeed, let ${\cal U}\subset \mathbb{S}^{d-2}$ be a small open domain
  in the unit sphere of dimension $d-2$. Then there exists a smooth $(d-1)\times(d-1)$ matrix-valued function $V(w)$,
  $w\in {\cal U}$, depending on ${\cal U}$, such that $V^{-1}(w)=V^t(w)$ and $V(w)w=\widetilde e_1=(1,0,...,0)\in \mathbb{R}^{d-1}$.
  Then we define the matrix $\Theta(\xi')$ for $\xi'/|\xi'|\in {\cal U}$ by
   $$\Theta(\xi')=\left(
\begin{array}{cc}
1&0\\
0&V(\xi'/|\xi'|)
\end{array}
\right).$$
It is easy to see that $\Theta(\xi')$ has the desired properties as long as $\xi'/|\xi'|\in {\cal U}$. Thus we can cover
$\mathbb{S}^{d-2}$ by a finite number of open sets ${\cal U}_k$, $k=1,...,K,$ so that to each ${\cal U}_k$ we can associate
a matrix-valued function $\Theta_k(\xi')$ having the desired properties for $\xi'/|\xi'|\in {\cal U}_k$. Then the identity
(\ref{eq:2.4}) remains valid with $\Theta(\xi')$ replaced by $\Theta_k(\xi')$ as long as $\xi'/|\xi'|\in {\cal U}_k$.
This implies the bounds (\ref{eq:2.5}) and (\ref{eq:2.7}) for $\xi'/|\xi'|\in {\cal U}_k$, $k=1,...,K,$ and hence for all 
$\xi'/|\xi'|\in \mathbb{S}^{d-2}$.
 \eproof

 \section{Some properties of the $h-\Psi$DOs} 

 We will first introduce the spaces of symbols which will play an important role in our analysis and will recall
 some basic properties of the $h-\Psi$DOs.
Given $k\in\mathbb{R}$, $\delta_1,\delta_2\ge 0$, we denote by 
$S_{\delta_1,\delta_2}^k$ the space of all functions $a\in C^\infty(T^*\Gamma)$, which may depend on the semiclassical parameter
$h$,  satisfying
$$\left|\partial_{x'}^\alpha\partial_{\xi'}^\beta a(x',\xi',h)\right|\le C_{\alpha,\beta}\langle\xi'\rangle^{k-\delta_1|\alpha|-\delta_2|\beta|}$$
for all multi-indices $\alpha$ and $\beta$, with constants $C_{\alpha,\beta}$ independent of $h$. 
More generally, given a function $\omega>0$ on $T^*\Gamma$, we denote by $S_{\delta_1,\delta_2}^k(\omega)$ the space of all functions $a\in C^\infty(T^*\Gamma)$, which may depend on the semiclassical parameter
$h$,  satisfying
$$\left|\partial_{x'}^\alpha\partial_{\xi'}^\beta a(x',\xi',h)\right|\le C_{\alpha,\beta}\omega^{k-\delta_1|\alpha|-\delta_2|\beta|}$$
for all multi-indices $\alpha$ and $\beta$, with constants $C_{\alpha,\beta}$ independent of $h$ and $\omega$. 
Thus $S_{\delta_1,\delta_2}^k=S_{\delta_1,\delta_2}^k(\langle\xi'\rangle)$.   
Given a matrix-valued symbol $a$, we will say that $a\in S_{\delta_1,\delta_2}^k$ if all entries of $a$
belong to $S_{\delta_1,\delta_2}^k$. 
Also, given $k\in\mathbb{R}$, $0\le\delta<1/2$, we denote by ${\cal S}_\delta^k$ the space of all
functions $a\in C^\infty(T^*\Gamma)$, which may depend on the semiclassical parameter
$h$, satisfying
$$\left|\partial_{x'}^\alpha\partial_{\xi'}^\beta a(x',\xi',h)\right|\le C_{\alpha,\beta}h^{-\delta(|\alpha|+|\beta|)}\langle\xi'\rangle^{k-|\beta|}$$
for all multi-indices $\alpha$ and $\beta$, with constants $C_{\alpha,\beta}$ independent of $h$.  
Again, given a matrix-valued symbol $a$, we will say that $a\in {\cal S}_\delta^k$ if all entries of $a$
belong to ${\cal S}_\delta^k$. 
The $h-\Psi$DO with a symbol $a$ is defined by
$$\left({\rm Op}_h(a)f\right)(x')=(2\pi h)^{-d+1}\int\int e^{-\frac{i}{h}\langle x'-y',\xi'\rangle}a(x',\xi',h)f(y')d\xi'dy'.$$
If $a\in S_{0,1}^k$, then the operator ${\rm Op}_h(a):H_h^k(\Gamma)\to L^2(\Gamma)$ is bounded uniformly in $h$,
where 
$$\left\|u\right\|_{H_h^k(\Gamma)}:= \left\|{\rm Op}_h(\langle\xi'\rangle^k)u\right\|_{L^2(\Gamma)}.$$ 
It is also well-known (e.g. see Section 7 of \cite{kn:DS}) that, if $a\in {\cal S}_\delta^0$, $0\le\delta<1/2$, then 
${\rm Op}_h(a):H_h^s(\Gamma)\to 
H_h^s(\Gamma)$ is bounded uniformly in $h$. More generally, we have the following 
\begin{prop} Let $h^{1/2-\epsilon}\le\theta\le 1$, $\ell\ge 0$, and let 
 $$a\in S_{1,1}^{-\ell}(\theta)+S_{0,1}^k\subset\theta^{-\ell}{\cal S}^k_{1/2-\epsilon}.$$
  Then we have
\begin{equation}\label{eq:3.1}
\left\|{\rm Op}_h(a)\right\|_{H_h^k(\Gamma)\to L^2(\Gamma)}\lesssim \theta^{-\ell}.
\end{equation}
\end{prop}
Let $\eta\in C^\infty(T^*\Gamma)$ be such that $\eta=1$ for $r_0\le C_0$, $\eta=0$ for $r_0\ge 2C_0$, where $C_0>0$ does not depend on $h$. 
Let $\rho$ denote either $\rho_s$ or $\rho_p$. 
It is easy to see (e.g. see Lemma 3.1 of \cite{kn:V1}) that taking $C_0$ big enough we can arrange 
$$C_1\theta^{1/2}\le |\rho|\le C_2,\quad {\rm Im}\,\rho\ge C_3|\theta||\rho|^{-1}\ge C_4|\theta|$$
 for $(x',\xi')\in{\rm supp}\,\eta$, and 
 $$|\rho|\ge {\rm Im}\,\rho\ge C_5|\xi'|$$
  for $(x',\xi')\in{\rm supp}\,(1-\eta)$ with some constants $C_j>0$. We will say that a function $a\in C^\infty(T^*\Gamma)$ belongs to
  $S_{\delta_1,\delta_2}^{k_1}(\omega_1)+S_{\delta_3,\delta_4}^{k_2}(\omega_2)$ if $\eta a\in S_{\delta_1,\delta_2}^{k_1}(\omega_1)$
  and $(1-\eta)a\in S_{\delta_3,\delta_4}^{k_2}(\omega_2)$. 
  It is shown in Section 3 of \cite{kn:V1} (see Lemma 3.2 of \cite{kn:V1}) that 
\begin{equation}\label{eq:3.2} 
\rho^k, |\rho|^k\in S_{2,2}^k(|\rho|)+S_{0,1}^k(|\rho|)\subset S_{1,1}^{-\widetilde k/2}(\theta)+S_{0,1}^k\subset \theta^{-\widetilde k/2}{\cal S}^{-N}_{1/2-\epsilon}+S_{0,1}^k
\subset\theta^{-\widetilde k/2}{\cal S}^k_{1/2-\epsilon}
\end{equation}
 as long as $\theta\ge h^{1/2-\epsilon}$, uniformly in $\theta$ and $h$, 
 where $\widetilde k=0$ if $k\ge 0$, $\widetilde k=-k$ if $k\le 0$ and $N\gg 1$ is arbitrary.

\section{A priori estimates}

In this section we will prove a priori estimates for the solution to the equation
\begin{equation}\label{eq:4.1}
\left\{
\begin{array}{l}
(h^2\Delta_{\lambda,\mu}+z^2n)u=hv\quad \mbox{in}\quad\Omega,\\
u=0\quad\mbox{on}\quad\Gamma.\\
\end{array}
\right.
\end{equation}

More precisely, we will prove the following

\begin{Theorem} \label{4.1} Let $\theta\ge h$ and $0<h\ll 1$. Let $u\in H^2(\Omega;\mathbb{C}^d)$ satisfy equation (\ref{eq:4.1}). 
Then the function $g=hB_{\lambda,\mu}u|_\Gamma$ satisfies the estimate
\begin{equation}\label{eq:4.2}
\|g\|_{L^2(\Gamma;\mathbb{C}^d)}\lesssim h^{1/2}\theta^{-1/2}\|v\|_{L^2(\Omega;\mathbb{C}^d)}.
\end{equation}
\end{Theorem}

{\it Proof.} We will first prove the following

\begin{lemma} \label{4.2} We have the estimate
\begin{equation}\label{eq:4.3}
\|u\|_{H_h^1(\Omega;\mathbb{C}^d)}\lesssim h\theta^{-1}\|v\|_{L^2(\Omega;\mathbb{C}^d)}.
\end{equation}
\end{lemma}

{\it Proof.} The analog of the Green formula for the elastic Laplacian applied to the solution $u$ of (\ref{eq:4.1}) takes the form
\begin{equation}\label{eq:4.4}
\left\langle -\Delta_{\lambda,\mu}u,u\right\rangle_{L^2(\Omega;\mathbb{C}^d)}=\int_\Omega E(u)
\end{equation}
where
$$E(u)=\lambda\sum_{j=1}^d\left|\frac{\partial u_j}{\partial x_j}\right|^2
+\frac{\mu}{2}\sum_{1\le i,j\le d}\left|\frac{\partial u_i}{\partial x_j}+\frac{\partial u_j}{\partial x_i}\right|^2$$
$$=(\lambda+2\mu)\sum_{j=1}^d\left|\frac{\partial u_j}{\partial x_j}\right|^2
+\frac{\mu}{2}\sum_{i\neq j}\left|\frac{\partial u_i}{\partial x_j}+\frac{\partial u_j}{\partial x_i}\right|^2$$
$$\ge C_1\sum_{1\le i,j\le d}\left|\frac{\partial u_i}{\partial x_j}+\frac{\partial u_j}{\partial x_i}\right|^2$$
with some constant $C_1>0$. On the other hand, since $u=0$ on $\Gamma$, by Korn's inequality we have
$$\int_\Omega\sum_{1\le i,j\le d}\left|\frac{\partial u_i}{\partial x_j}\right|^2
\le C_2\int_\Omega\sum_{1\le i,j\le d}\left|\frac{\partial u_i}{\partial x_j}+\frac{\partial u_j}{\partial x_i}\right|^2$$
with some constant $C_2>0$. Combining the above inequalities with (\ref{eq:4.4}) we obtain the coercive estimate
\begin{equation}\label{eq:4.5}
\left\langle -\Delta_{\lambda,\mu}u,u\right\rangle_{L^2(\Omega;\mathbb{C}^d)}\ge 
C\int_\Omega\sum_{1\le i,j\le d}\left|\frac{\partial u_i}{\partial x_j}\right|^2
\end{equation}
with some constant $C>0$. The Green formula (\ref{eq:4.4}) also gives the identity
$${\rm Im}(z^2)\left\|n^{1/2}u\right\|^2_{L^2(\Omega;\mathbb{C}^d)}={\rm Im}\,\left\langle hv,u\right\rangle_{L^2(\Omega;\mathbb{C}^d)},$$
which implies
\begin{equation}\label{eq:4.6}
\|u\|_{L^2(\Omega;\mathbb{C}^d)}\lesssim h\theta^{-1}\|v\|_{L^2(\Omega;\mathbb{C}^d)}.
\end{equation}
On the other hand, we have
$$\left\langle -h^2\Delta_{\lambda,\mu}u,u\right\rangle_{L^2(\Omega;\mathbb{C}^d)}={\rm Re}(z^2)\left\langle nu,u\right\rangle_{L^2(\Omega;\mathbb{C}^d)}-{\rm Re}\,\left\langle hv,u\right\rangle_{L^2(\Omega;\mathbb{C}^d)}$$
$$\lesssim \|u\|^2_{L^2(\Omega;\mathbb{C}^d)}+h^2\|v\|^2_{L^2(\Omega;\mathbb{C}^d)}$$
which combinned with (\ref{eq:4.5}) leads to the estimate 
\begin{equation}\label{eq:4.7}
\int_\Omega\sum_{1\le i,j\le d}h^2\left|\frac{\partial u_i}{\partial x_j}\right|^2\lesssim
\|u\|^2_{L^2(\Omega;\mathbb{C}^d)}+ h^2\|v\|^2_{L^2(\Omega;\mathbb{C}^d)}.
\end{equation}
Clearly, (\ref{eq:4.3}) follows from (\ref{eq:4.6}) and (\ref{eq:4.7}).
\eproof

Let ${\cal V}\subset\mathbb{R}^d$ be a small open domain such that ${\cal V}^0:={\cal V}\cap\Gamma\neq\emptyset$. 
Let $(x_1,x')\in {\cal V}^+:={\cal V}\cap\Omega$, $0<x_1\ll 1$, $x'=(x_2,...,x_d)\in{\cal V}^0$, be the local normal geodesic coordinates near the boundary. Recall (e.g. see Section 2 of \cite{kn:V3}) that the Euclidean gradient $\nabla$ can be written in the coordinates $x=(x_1,x')$ as 
$$\nabla=\gamma(x)\nabla_x=\nu(x')\frac{\partial}{\partial x_1}+\sum_{k=2}^d\gamma(x)e_k\frac{\partial}{\partial x_k},$$ 
where $\gamma$ is a smooth matrix-valued function such that $\gamma(x)e_1=\nu(x')$,
$\gamma(x)e_k$ satisfy 
 \begin{equation}\label{eq:4.8}
 \langle \nu(x'),\gamma(x)e_k\rangle =0,\quad k=2,...,d.
\end{equation}
Let $\xi=(\xi_1,\xi')$ be the dual variable of $x=(x_1,x')$. 
Then the symbol of the operator $-i\nabla$
in the coordinates $(x,\xi$) takes the form $\xi_1\nu(x')+\beta(x,\xi')$, where 
$$\beta(x,\xi')=\sum_{k=2}^d\xi_k\gamma(x)e_k.$$
Note that (\ref{eq:4.8}) implies the identity
\begin{equation}\label{eq:4.9}
 \langle \nu(x'),\beta(x,\xi')\rangle =0\quad \mbox{for all}\quad (x,\xi').
\end{equation}
Thus we get that the principal symbol of $-\Delta$ is equal to 
$\xi_1^2+r(x,\xi')$, where $r=\langle \beta,\beta\rangle$. 
Therefore, the principal symbol of the positive Laplace-Beltrami operator on $\Gamma$ is equal to
$r_0(x',\xi')=r(0,x',\xi')=\langle \beta_0,\beta_0\rangle$, where $\beta_0=\beta|_{x_1=0}$. Clearly, there exist constants
$C_1,C_2>0$ such that $C_1|\xi'|^2\le r_0\le C_2|\xi'|^2$. 

Let ${\cal V}_1\subset{\cal V}$ be a small open domain such that ${\cal V}_1^0:={\cal V_1}\cap\Gamma\neq\emptyset$. 
Choose a function $\psi\in C_0^\infty({\cal V})$, $0\le\psi\le 1$, such that $\psi=1$ on ${\cal V}_1$. 
 Then the function $u^\flat:=\psi u$ satisfies the equation 
\begin{equation}\label{eq:4.10}
\left\{
\begin{array}{l}
(h^2\Delta_{\lambda,\mu}+z^2n)u^\flat=hv^\flat\quad \mbox{in}\quad\Omega,\\
u^\flat=0\quad\mbox{on}\quad\Gamma,\\
\end{array}
\right.
\end{equation}
where $v^\flat=\psi v+h[\Delta_{\lambda,\mu},\psi]u$ satisfies
\begin{equation}\label{eq:4.11}
\|v^\flat\|_{L^2(\Omega;\mathbb{C}^d)}\lesssim\|v\|_{L^2(\Omega;\mathbb{C}^d)}+\|u\|_{H^1_h(\Omega;\mathbb{C}^d)}.
\end{equation}
 We will now write the elastic Laplacian in the coordinates $x=(x_1,x')$. To this end, we will write the principal symbol of 
 $-\Delta_{\lambda,\mu}$ in the coordinates $(x,\xi)$. We have
 $$P(x,\xi)=\mu(\gamma\xi)^2I_d+(\lambda+\mu)(\gamma\xi)\otimes(\gamma\xi)$$
 $$=\mu(\xi_1^2+r(x,\xi'))I_d+(\lambda+\mu)\gamma(\xi\otimes\xi)\gamma^t$$
 $$=\xi_1^2Q_0(x)+\xi_1Q_1(x,\xi')+Q_2(x,\xi'),$$
 where
 $$Q_0=c_s\Pi_s(e_1)+c_p\Pi_p(e_1),$$
$$Q_1=(\lambda+\mu)\gamma\left(e_1\otimes\xi'+\xi'\otimes e_1\right)\gamma^t,$$
 $$Q_2=\mu r(x,\xi')I_d+(\lambda+\mu)\gamma(\xi'\otimes\xi')\gamma^t$$
are symmetric matrices. Denote ${\cal D}_{x_j}=-ih\partial_{x_j}$. We can write
\begin{equation}\label{eq:4.12}
-h^2\Delta_{\lambda,\mu}=Q_0(x){\cal D}_{x_1}^2+{\cal Q_1}{\cal D}_{x_1}+{\cal Q_2}+h{\cal R}(x,{\cal D}_{x}),
\end{equation}
where ${\cal R}$ is a first-order matrix-valued differential operator, and 
$${\cal Q_j}=\frac{1}{2}\left(Q_j(x,{\cal D}_{x'})+Q_j(x,{\cal D}_{x'})^*\right)=Q_j(x,{\cal D}_{x'})
+h{\cal R}_{j-1}(x,{\cal D}_{x'}),\quad j=1,2,$$
are self-adjoint operators on $L^2(\Gamma;\mathbb{C}^d)$. Here $Q^*$ denotes the adjoint of $Q$ with respect to the scalar product,
$\langle\cdot,\cdot\rangle_0$, 
in $L^2(\Gamma;\mathbb{C}^d)$, and ${\cal R}_{j-1}$ is a $j-1$-order matrix-valued differential operator. Introduce the function
$$F(x_1)=\left\langle Q_0(x_1,\cdot){\cal D}_{x_1}u^\flat,{\cal D}_{x_1}u^\flat\right\rangle_0
-\left\langle {\cal Q_2}(x_1,\cdot,{\cal D}_{x'})u^\flat,u^\flat\right\rangle_0+{\rm Re}(z^2)\left\langle n(x_1,\cdot)u^\flat,u^\flat\right\rangle_0.$$
Clearly,
\begin{equation}\label{eq:4.13}
F(0)=\left\langle Q_0(0,\cdot){\cal D}_{x_1}u^\flat|_{x_1=0},{\cal D}_{x_1}u^\flat|_{x_1=0}\right\rangle_0\ge 
C\left\|{\cal D}_{x_1}u^\flat|_{x_1=0}\right\|_0^2
\end{equation}
with some constant $C>0$, where $\|\cdot\|_0$ denotes the norm in $L^2(\Gamma;\mathbb{C}^d)$. On the other hand,
\begin{equation}\label{eq:4.14}
F(0)=-\int_0^{\delta}F'(x_1)dx_1
\end{equation}
for some constant $\delta>0$, 
where $F'$ denotes the first derivative with respect to $x_1$. We will now use (\ref{eq:4.14}) to bound $F(0)$ from above. 
To this end we will compute $F'(x_1)$ using that $u^\flat$ satisfies (\ref{eq:4.10}) together with (\ref{eq:4.12}).
We have
$$F'(x_1)=-2{\rm Re}\,\left\langle (Q_0{\cal D}_{x_1}^2+{\cal Q_2}-{\rm Re}(z^2)n)u^\flat,\partial_{x_1}u^\flat\right\rangle_0$$
$$+\left\langle Q'_0{\cal D}_{x_1}u^\flat,{\cal D}_{x_1}u^\flat\right\rangle_0
-\left\langle ({\cal Q_2}'-{\rm Re}(z^2)n')u^\flat,u^\flat\right\rangle_0$$
$$=2h^{-1}{\rm Im}\,\left\langle (h^2\Delta_{\lambda,\mu}+{\rm Re}(z^2)n)u^\flat,{\cal D}_{x_1}u^\flat\right\rangle_0+
2h^{-1}{\rm Im}\,\left\langle (Q_1{\cal D}_{x_1}+h{\cal R})u^\flat,{\cal D}_{x_1}u^\flat\right\rangle_0$$
$$+\left\langle Q'_0{\cal D}_{x_1}u^\flat,{\cal D}_{x_1}u^\flat\right\rangle_0
-\left\langle ({\cal Q_2}'-{\rm Re}(z^2)n')u^\flat,u^\flat\right\rangle_0$$
$$=2{\rm Im}\,\left\langle (v^\flat-ih^{-1}{\rm Im}(z^2)nu^\flat),{\cal D}_{x_1}u^\flat\right\rangle_0+
2{\rm Im}\,\left\langle {\cal R}u^\flat,{\cal D}_{x_1}u^\flat\right\rangle_0$$
$$+\left\langle Q'_0{\cal D}_{x_1}u^\flat,{\cal D}_{x_1}u^\flat\right\rangle_0
-\left\langle ({\cal Q_2}'-{\rm Re}(z^2)n')u^\flat,u^\flat\right\rangle_0.$$
Hence
$$|F'(x_1)|\lesssim h\theta^{-1}\|v^\flat\|_0^2+\theta h^{-1}\sum_{\ell=0}^1\|{\cal D}_{x_1}^\ell u^\flat\|_0^2+\sum_{|\alpha|\le 1}\|{\cal D}_{x}^\alpha  u^\flat\|_0^2.$$
Using this estimate together with (\ref{eq:4.11}), (\ref{eq:4.14}) and Lemma 4.2 we obtain
\begin{equation}\label{eq:4.15}
F(0)\le \int_0^{2\delta}|F'(x_1)|dx_1\lesssim h\theta^{-1}\|v\|_{L^2(\Omega;\mathbb{C}^d)}^2
+(1+\theta h^{-1})\|u\|_{H_h^1(\Omega;\mathbb{C}^d)}^2$$ 
$$\lesssim (h\theta^{-1}+h^2\theta^{-2})\|v\|_{L^2(\Omega;\mathbb{C}^d)}^2\lesssim h\theta^{-1}\|v\|_{L^2(\Omega;\mathbb{C}^d)}^2.
\end{equation}
Observe now that
 $${\cal D}_{x_1}u^\flat|_{x_1=0}=\psi_0{\cal D}_{x_1}u|_{x_1=0},\quad {\cal D}_{x'}u|_{x_1=0}=0,$$
 where $\psi_0=\psi|_{x_1=0}$ is supported in ${\cal V}^0$ and such that $\psi_0=1$ on ${\cal V}_1^0$. 
 Therefore, 
  by (\ref{eq:4.13}) and (\ref{eq:4.15}),
$$
\left\|\psi_0{\cal D}_{x_1}u|_{x_1=0}\right\|_0\lesssim h^{1/2}\theta^{-1/2}\|v\|_{L^2(\Omega;\mathbb{C}^d)},
$$
which clearly implies 
\begin{equation}\label{eq:4.16}
\left\|\psi_0g\right\|_0\lesssim h^{1/2}\theta^{-1/2}\|v\|_{L^2(\Omega;\mathbb{C}^d)}.
\end{equation}
Since $\Gamma$ is compact, there exist a finite number of smooth functions $\psi_i$, $0\le\psi_i\le 1$, $i=1,...,I,$ such that 
$1=\sum_{i=1}^I\psi_i$ and (\ref{eq:4.16}) holds with $\psi_0$ replaced by each $\psi_i$. Therefore, 
 the estimate (\ref{eq:4.2}) is obtained by summing up all such estimates (\ref{eq:4.16}).
\eproof

\section{Parametrix construction} 

We keep the notations from the previous sections and will suppose that $\theta\ge h^{2/5-\epsilon}$,
$0<\epsilon\ll 1$.  
It suffices to build the parametrix locally since the global parametrix can be obtained by using a suitable partition of the
unity and summing up the corresponding local parametrices. 
Let the function $\phi_0\in C_0^\infty(\mathbb{R})$ 
be such that $\phi_0(\sigma)=1$  for $|\sigma|\le 1$, $\phi_0(\sigma)=0$  for $|\sigma|\ge 2$.
Let $(x_1,x')\in {\cal V}^+$ be the local normal geodesic coordinates near the boundary. 
Take a function $\chi\in C^\infty(T^*\Gamma)$, $0\le\chi\le 1$, such that $\pi_{x'}({\rm supp}\,\chi)\subset {\cal V}^0$, where 
$\pi_{x'}:T^*\Gamma\to\Gamma$ denotes the projection $(x',\xi')\to x'$. Moreover, we require that either $\chi$ is of compact support or
 $\chi\in S_{0,1}^0$ with ${\rm supp}\,\chi\subset{\rm supp}(1-\eta)$. When $\chi$ is of compact support we require that ${\rm supp}\,\chi$
 has common points with at most one glancing region. 
Let $f\in H^3(\Gamma;\mathbb{C}^d)$.  
We will be looking for a parametrix of the solution to 
equation (\ref{eq:1.1}) in the form
$$\widetilde u=(2\pi h)^{-d+1}\int\int e^{\frac{i}{h}(\langle y',\xi'\rangle+\varphi_s(x,\xi',z))}\Psi(x,\xi')A_s(x,\xi',z)\chi(x',\xi')f(y')d\xi'dy'$$
$$+(2\pi h)^{-d+1}\int\int e^{\frac{i}{h}(\langle y',\xi'\rangle+\varphi_p(x,\xi',z))}\Psi(x,\xi')A_p(x,\xi',z)\chi(x',\xi')f(y')d\xi'dy',$$
where 
$$\Psi=\phi_0(x_1\langle\xi'\rangle^\varepsilon/\delta)\phi_0\left(x_1/\left|\rho_s\right|^3\delta\right)\phi_0\left(x_1/\left|\rho_p\right|^3\delta\right),\quad 0<\varepsilon\ll 1,$$
$0<\delta\ll 1$ being a parameter independent of $h$ and $\theta$ to be fixed in Lemma 5.1.  
 We require that $\widetilde u$ satisfies the boundary condition
$\widetilde u={\rm Op}_h(\chi)f$ on $x_1=0$. The phase functions are of the form
$$\varphi_s=\sum_{k=0}^{N-1}x_1^k\varphi_{s,k},\quad \varphi_{s,0}=-\langle x',\xi'\rangle,\,\varphi_{s,1}=\rho_s,$$
$$\varphi_p=\sum_{k=0}^{N-1}x_1^k\varphi_{p,k},\quad \varphi_{p,0}=-\langle x',\xi'\rangle,\,\varphi_{p,1}=\rho_p,$$
$N\gg 1$ being an arbitrary integer, and satisfy the eikonal equations mod ${\cal O}(x_1^N)$:
\begin{equation}\label{eq:5.1}
\left\{
\begin{array}{l}
c_s(x)(\gamma\nabla_x\varphi_s)^2 -z^2n(x)=x_1^N\Phi_s,\\
c_p(x)(\gamma\nabla_x\varphi_p)^2-z^2n(x)=x_1^N\Phi_p,
\end{array}
\right.
\end{equation}
where $\Phi_s$, $\Phi_p$ are smooth functions up to the boundary $x_1=0$. 
  One can solve the eikonal equations above in the same way as in \cite{kn:V1}.
The functions $\varphi_{s,k}$, $\varphi_{p,k}$, $k\ge 2$, are determined uniquely,
independent of $x_1$, and have the following properties (see Section 4 of \cite{kn:V1}). 
 
\begin{lemma} For  $0\le x_1\le 2\delta\min\{1,|\rho_s|^3\}$ with $\delta>0$ small enough, we have
\begin{equation}\label{eq:5.2}
\varphi_{s,k}\in S_{2,2}^{4-3k}\left(\left|\rho_s\right|\right)+S_{0,1}^1,\quad k\ge 1,
\end{equation}
\begin{equation}\label{eq:5.3}
\partial_{x_1}^k\Phi_s\in S_{2,2}^{2-3N-3k}\left(\left|\rho_s\right|\right)+S_{0,1}^2,\quad k\ge 0,
\end{equation}
\begin{equation}\label{eq:5.4}
{\rm Im}\,\varphi_s\ge x_1{\rm Im}\,\rho_s/2,
\end{equation} 
\begin{equation}\label{eq:5.5}
\left|\partial_{x_1}\varphi_s\right|\ge \left|\rho_s\right|/2,
\end{equation}
and similarly for $\varphi_p$. 
\end{lemma}

Set $\widetilde\varphi_s=\varphi_s-\varphi_{s,0}$, $\widetilde\varphi_p=\varphi_p-\varphi_{p,0}$. The next lemma is proved in Section 4 of
\cite{kn:V3}.

\begin{lemma} There exists a constant $C>0$ such that we have the estimates
\begin{equation}\label{eq:5.6} 
\left|\partial_{x'}^\alpha\partial_{\xi'}^\beta\left(e^{i\widetilde\varphi_s/h}\right)\right|\le 
\left\{
\begin{array}{l}
C_{\alpha,\beta}\theta^{-|\alpha|-|\beta|}e^{-Cx_1\theta/h}
\quad\mbox{on}\quad{\rm supp}\,\eta,\\
C_{\alpha,\beta}|\xi'|^{-|\beta|}e^{-Cx_1|\xi'|/h}
\quad\mbox{on}\quad{\rm supp}(1-\eta),
\end{array}
\right.
\end{equation}
for $0\le x_1\le 2\delta\min\{1,|\rho_s|^3\}$ and all multi-indices $\alpha$ and $\beta$ with constants $C_{\alpha,\beta}>0$ independent of $x_1$, $\theta$, $z$ and $h$. Similar bounds hold for $\widetilde\varphi_p$ as well. 
\end{lemma}

 The amplitudes $A_s$ and $A_p$ are matrix-valued functions which will be chosen so that on supp$\,\chi$ we have 
 \begin{equation}\label{eq:5.7}
 A_s+A_p=I_d\quad\mbox{on}\quad x_1=0,
 \end{equation}
  and
 \begin{equation}\label{eq:5.8} 
 \Pi_p(\gamma\nabla_x\varphi_s)A_s=0,\quad \Pi_s(\gamma\nabla_x\varphi_p)A_p=0.
 \end{equation}
 If ${\cal V}^0$ is small enough, there exists 
a matrix-valued function $\Lambda(x')\in C^\infty({\cal V}^0)$ such that $\Lambda(x')\nu(x')=e_1$ and $\Lambda^t=\Lambda^{-1}$
in ${\cal V}^0$. Set
$$U(\xi)=\Lambda^{-1}U_0(\Lambda\xi)\Lambda,$$
where $U_0$ is the matrix introduced in Section 2. It follows from (\ref{eq:2.2}) that $U(\xi)$ 
is invertible if $\langle \xi,\nu\rangle\neq 0$. Clearly, $U(\nu)=I_d$. Moreover, by Lemma 2.1 we have $U(\xi)\xi=\xi^2\nu$. Therefore,
 \begin{equation}\label{eq:5.9} 
 U(\xi)\Pi_p(\xi)U^t(\xi)=U(\xi)(\xi\otimes\xi)U^t(\xi)
 =(U(\xi)\xi)\otimes(U(\xi)\xi)=\xi^4\nu\otimes\nu=\xi^4\Pi_p(\nu).
 \end{equation}
 Moreover, we have
 $$Z(\xi):=U(\xi)U^t(\xi)=\Lambda^{-1}Z_0(\Lambda\xi)\Lambda,\quad \Pi_p(\nu)=\Lambda^{-1}\Pi_p(e_1)\Lambda.$$
 Hence the matrices $Z(\xi)$ and $\Pi_p(\nu)$ commute. Set 
 $$A_s=U^t(\gamma\nabla_x\varphi_s)\Pi_s(\nu)T,\quad A_p=U^t(\gamma\nabla_x\varphi_p)\Pi_p(\nu)T,$$
 where $T$ is a matrix-valued function independent of $x_1$ to be defined below in such a way that (\ref{eq:5.7}) holds. 
 Let us see that $A_s$ and $A_p$ satisfy (\ref{eq:5.8}).  In view of (\ref{eq:5.5}) we have 
 $$\langle\gamma\nabla_x\varphi_s,\nu\rangle=\partial_{x_1}\varphi_s\neq 0,\quad  
 \langle\gamma\nabla_x\varphi_p,\nu\rangle=
 \partial_{x_1}\varphi_p\neq 0.$$
 Hence the matrices $U^t(\gamma\nabla_x\varphi_s)$ and $U^t(\gamma\nabla_x\varphi_p)$ are invertible, and by (\ref{eq:5.9}), we have 
 $$\Pi_p(\gamma\nabla_x\varphi_s)A_s=\Pi_p(\gamma\nabla_x\varphi_s)U^t(\gamma\nabla_x\varphi_s)\Pi_s(\nu)T$$ 
 $$=(\gamma\nabla_x\varphi_s)^4U^{-1}(\gamma\nabla_x\varphi_s)\Pi_p(\nu)\Pi_s(\nu)T=0,$$
 $$A_p=U^{-1}(\gamma\nabla_x\varphi_p)Z(\gamma\nabla_x\varphi_p)\Pi_p(\nu)T=U^{-1}(\gamma\nabla_x\varphi_p)\Pi_p(\nu) 
 Z(\gamma\nabla_x\varphi_p)T$$
 $$=(\gamma\nabla_x\varphi_p)^{-4}\Pi_p(\gamma\nabla_x\varphi_p)Z(\gamma\nabla_x\varphi_p)T,$$
 which imply (\ref{eq:5.8}). We will now find the matrix $T$ so that  
 $WT=I_d$ with
 $$W=U^t(\gamma\nabla_x\varphi_s|_{x_1=0})\Pi_s(\nu)+U^t(\gamma\nabla_x\varphi_p|_{x_1=0})\Pi_p(\nu).$$
 Observe that
$$\gamma\nabla_x\varphi_s|_{x_1=0}=\rho_s\nu-\beta_0,\quad\gamma\nabla_x\varphi_p|_{x_1=0}=\rho_p\nu-\beta_0.$$
 Therefore,
 $$W=\rho_s\Pi_s(\nu)+\rho_p\Pi_p(\nu)-U^t(\beta_0).$$
 We will derive from Lemma 2.2 the following
 
  \begin{lemma} The matrix $W$ is invertible with an inverse $T=W^{-1}$ satisfying the bounds
  \begin{equation}\label{eq:5.10}
  \left\|T\right\|\lesssim 
  \left\{
\begin{array}{l}
\theta^{-\ell}\quad\mbox{on}\quad{\rm supp}\,\eta,\\
\sqrt{r_0+1}\quad\mbox{on}\quad{\rm supp}(1-\eta),
\end{array}
\right.
   \end{equation}
   where $\ell=0$ if $d=2$, $\ell=1/2$ if $d\ge 3$.
   More generally, we have
    \begin{equation}\label{eq:5.11}
    T\in S_{1,1}^{-\ell}(\theta)+S_{0,1}^1.
     \end{equation}
   \end{lemma}
   
   {\it Proof.} Set 
 $\zeta(x',\xi')=-\Lambda(x')\beta_0(x',\xi')$. Clearly, $\zeta^2=\beta_0^2=r_0$. 
 In view of (\ref{eq:4.9}) we have $\langle \zeta,e_1\rangle=0$ and hence
 $\zeta=(0,\zeta_2,...,\zeta_d)$. 
 Then the matrix $W_0:=\Lambda W\Lambda^{-1}$ can be written in the form
 $$W_0=U_0^t(\zeta)+\rho_s\Pi_s(e_1)+\rho_p\Pi_p(e_1)$$
 $$=U_0^t(\rho_se_1+\zeta)+(\rho_p-\rho_s)e_1\otimes e_1.$$
 By (\ref{eq:2.2}) we get
   \begin{equation}\label{eq:5.12}
   {\rm det}W_0=(r_0+\rho_s\rho_p)\rho_s^{d-2}.
     \end{equation}
     We need now the following
     
      \begin{lemma} There exists a constant $C>0$ such that 
     \begin{equation}\label{eq:5.13}
   |r_0+\rho_s\rho_p|\ge C.
     \end{equation}
      More generally, we have
    \begin{equation}\label{eq:5.14}
    (r_0+\rho_s\rho_p)^{-1}\in S_{1,1}^{0}(\theta)+S_{0,1}^0.
     \end{equation}
      \end{lemma}
      
      {\it Proof.} Recall that $\rho_s^2=-r_0+z^2k_s$, $\rho_p^2=-r_0+z^2k_p$ with some functions 
     $k_s,k_p\in C^\infty(\Gamma)$, $k_s>k_p>0$. Then we have the identity
     \begin{equation}\label{eq:5.15}
     r_0+\rho_s\rho_p=r_0+\rho_s^2+\frac{\rho_s(\rho_p^2-\rho_s^2)}{\rho_s+\rho_p}=z^2\frac{k_p\rho_s+k_s\rho_p}{\rho_s+\rho_p}.
      \end{equation}
     Hence 
      \begin{equation}\label{eq:5.16}
   |r_0+\rho_s\rho_p|\ge \frac{k_p{\rm Im}\,\rho_s+k_s{\rm Im}\,\rho_p}{|\rho_s|+|\rho_p|}\ge C_1\frac{{\rm Im}\,\rho_s
   +{\rm Im}\,\rho_p}{\sqrt{r_0+1}}
     \end{equation}
     with some constant $C_1>0$. On the other hand, there is a constant $C_2>0$ such that ${\rm Im}\,\rho_s\ge C_2\sqrt{r_0+1}$
     on supp$(1-\eta)$, ${\rm Im}\,\rho_s\ge C_2\theta$ on supp$\eta$, and similarly for $\rho_p$ (see Section 3). 
     Therefore, (\ref{eq:5.13}) follows from (\ref{eq:5.16}) when either $(x',\xi')\in{\rm supp}(1-\eta)$ or 
     $(x',\xi')\in{\rm supp}\,\eta$ and $\theta\ge\theta_0>0$. Thus, it remains to prove (\ref{eq:5.13}) when 
     $(x',\xi')\in{\rm supp}\,\eta$ and $\theta\ll 1$. In this case we have that $|\rho_s|$
     and $|\rho_p|$ are uniformly bounded from above by a constant and $z^2=1+{\cal O}(\theta)$. We will make use of the identity
      \begin{equation}\label{eq:5.17}
     (k_p\rho_s-k_s\rho_p)(k_p\rho_s+k_s\rho_p)=(k_s-k_p)\left((k_s+k_p)r_0-z^2k_sk_p\right).
     \end{equation}
     If 
     $$\left|(k_s+k_p)r_0-z^2k_sk_p\right|\ge\varepsilon>0,$$
     it follows from (\ref{eq:5.17}) that 
     $$|k_p\rho_s+k_s\rho_p|\ge C_3\varepsilon$$
     with some constant $C_3>0$. Thus in this case (\ref{eq:5.13}) follows from (\ref{eq:5.15}). Let now
     $$\left|(k_s+k_p)r_0-z^2k_sk_p\right|\le\varepsilon$$
     with $0<\varepsilon\ll 1$. Then 
     $$\rho_s=\sqrt{\frac{k_sr_0}{k_p}}+{\cal O}(\varepsilon),\quad \rho_p=\sqrt{\frac{k_pr_0}{k_s}}+{\cal O}(\varepsilon),\quad 
      r_0=\frac{k_sk_p}{k_s+k_p}+{\cal O}(\varepsilon)+{\cal O}(\theta).$$
     Clearly, there exists a constant $C_4>0$ such that $|k_p\rho_s+k_s\rho_p|\ge C_4$, provided $\varepsilon$
     and $\theta$ are taken small enough, which again implies (\ref{eq:5.13}). 
     
     To prove (\ref{eq:5.14}) note first that, in view of (\ref{eq:3.2}), we have 
     \begin{equation}\label{eq:5.18}
    \rho_s+\rho_p,\,k_p\rho_s+k_s\rho_p\in S_{1,1}^{0}(\theta)+S_{0,1}^1.
     \end{equation}
     Therefore, in view of (\ref{eq:5.15}), to prove (\ref{eq:5.14}) it suffices to show that
     \begin{equation}\label{eq:5.19}
    (k_p\rho_s+k_s\rho_p)^{-1}\in S_{1,1}^{0}(\theta)+S_{0,1}^{-1}.
     \end{equation}
     In other words, we must show that given any multi-indices $\alpha$ and $\beta$ we have the estimates
     \begin{equation}\label{eq:5.20} 
\left|\partial_{x'}^\alpha\partial_{\xi'}^\beta\left((k_p\rho_s+k_s\rho_p)^{-1}\right)\right|\le 
\left\{
\begin{array}{l}
C_{\alpha,\beta}\theta^{-|\alpha|-|\beta|}
\quad\mbox{on}\quad{\rm supp}\,\eta,\\
C_{\alpha,\beta}|\xi'|^{-1-|\beta|}
\quad\mbox{on}\quad{\rm supp}(1-\eta).
\end{array}
\right.
\end{equation}
Clearly, for $\alpha=\beta=0$ the bounds in (\ref{eq:5.20}) follow from the analysis above. To prove them for all
$\alpha$ and $\beta$, we will proceed by induction in $|\alpha|+|\beta|$. Suppose that (\ref{eq:5.20}) holds for all 
$\alpha$ and $\beta$ such that $|\alpha|+|\beta|\le K-1$, $K\ge 1$. Let us see that (\ref{eq:5.20}) holds for all 
$\alpha$ and $\beta$ such that $|\alpha|+|\beta|=K$. To this end, we will use the identity
$$0=\partial_{x'}^\alpha\partial_{\xi'}^\beta\left((k_p\rho_s+k_s\rho_p)(k_p\rho_s+k_s\rho_p)^{-1}\right)$$
$$=(k_p\rho_s+k_s\rho_p)\partial_{x'}^\alpha\partial_{\xi'}^\beta\left((k_p\rho_s+k_s\rho_p)^{-1}\right)$$
 $$+\sum_{|\alpha'|+|\beta'|\le K-1}\partial_{x'}^{\alpha'}\partial_{\xi'}^{\beta'}\left((k_p\rho_s+k_s\rho_p)^{-1}\right)
 \partial_{x'}^{\alpha-\alpha'}\partial_{\xi'}^{\beta-\beta'}\left(k_p\rho_s+k_s\rho_p\right). $$
 Thus, in view of (\ref{eq:5.18}), we obtain
 $$|k_p\rho_s+k_s\rho_p|\left|\partial_{x'}^\alpha\partial_{\xi'}^\beta\left((k_p\rho_s+k_s\rho_p)^{-1}\right)\right|$$
 $$\le \sum_{|\alpha'|+|\beta'|\le K-1}\left|\partial_{x'}^{\alpha'}\partial_{\xi'}^{\beta'}\left((k_p\rho_s+k_s\rho_p)^{-1}\right)\right|
 \left|\partial_{x'}^{\alpha-\alpha'}\partial_{\xi'}^{\beta-\beta'}\left(k_p\rho_s+k_s\rho_p\right)\right|$$
 $$\lesssim \left\{
\begin{array}{l}
\theta^{-|\alpha|-|\beta|}
\quad\mbox{on}\quad{\rm supp}\,\eta,\\
|\xi'|^{-|\beta|}
\quad\mbox{on}\quad{\rm supp}(1-\eta).
\end{array}
\right.$$
Since $|\xi'|+1\lesssim|k_p\rho_s+k_s\rho_p|$, we conclude from the above bounds that (\ref{eq:5.20}) holds for all 
$\alpha$ and $\beta$ such that $|\alpha|+|\beta|=K$, as desired.
  \eproof
 
 It follows from (\ref{eq:5.12}) and (\ref{eq:5.13}) that the matrix $W_0$ is invertible, and hence so is $W$. Moreover, by 
 (\ref{eq:2.3}) its inverse satisfies the bound
 $$\|T\|\lesssim\left\|W_0^{-1}\right\|\lesssim |\zeta|+|\rho_s|+|\rho_p|+(d-2)|\rho_s|^{-1}$$
 $$\lesssim \left\{
\begin{array}{l}
1+(d-2)\theta^{-1/2}
\quad\mbox{on}\quad{\rm supp}\,\eta,\\
\sqrt{r_0+1}
\quad\mbox{on}\quad{\rm supp}(1-\eta),
\end{array}
\right.$$
 which implies (\ref{eq:5.10}). To prove (\ref{eq:5.11}) we need to show that the estimates
     \begin{equation}\label{eq:5.21} 
\left\|\partial_{x'}^\alpha\partial_{\xi'}^\beta T\right\|\le 
\left\{
\begin{array}{l}
C_{\alpha,\beta}\theta^{-\ell-|\alpha|-|\beta|}
\quad\mbox{on}\quad{\rm supp}\,\eta,\\
C_{\alpha,\beta}|\xi'|^{1-|\beta|}
\quad\mbox{on}\quad{\rm supp}(1-\eta).
\end{array}
\right.
\end{equation}
 hold for all multi-indices $\alpha$ and $\beta$. Note that in view of (\ref{eq:3.2}) we have
   \begin{equation}\label{eq:5.22}
    W\in S_{1,1}^{0}(\theta)+S_{0,1}^1.
     \end{equation}
     Now (\ref{eq:5.21}) can be derived from (\ref{eq:5.10}) and (\ref{eq:5.22}) by induction in
     $|\alpha|+|\beta|$ in the same way as above.
  \eproof
 
  To get a parametrix for the elastic DN map we need the following
  
  \begin{lemma} There exist matrix-valued functions $m_d,q\in C^\infty(T^*\Gamma)$ such that
  \begin{equation}\label{eq:5.23}
  -ihB_{\lambda,\mu}\widetilde u|_{x_1=0}={\rm Op}_h(m_d\chi+hq\chi)f. 
   \end{equation}
   \end{lemma}
  
  {\it Proof.} Given a scalar-valued function $\varphi$ and a vector-valued function $a$, we have the identity
  $$-ihe^{-i\varphi/h}B_{\lambda,\mu}\left(e^{i\varphi/h}a\right)
=\lambda\langle \gamma\nabla_x\varphi,a\rangle\nu+\mu\langle\nu,a\rangle\gamma\nabla_x\varphi+
\mu\langle\nu,\gamma\nabla_x\varphi\rangle a$$
$$-ih\lambda\langle \gamma\nabla_x,a\rangle\nu-ih\mu\langle\nu,\gamma\nabla_xa\rangle-ih\mu\langle\nu,\gamma\nabla_x\rangle a.$$
Set $a_s=A_s\chi f$, $a_p=A_p\chi f$, $a_s^0=A_s^0\chi f$, $a_p^0=A_p^0\chi f$, where $A_s^0=A_s|_{x_1=0}$, $A_p^0=A_p|_{x_1=0}$ satisfy $A_s^0+A_p^0=I$. 
Set also $a_s^1=A_s^1\chi f$, $a_p^1=A_p^1\chi f$, where $A_s^1=\partial_{x_1}A_s|_{x_1=0}$, $A_p^1=\partial_{x_1}A_p|_{x_1=0}$. 
Applying the above identity to $\varphi_s$, $a_s$ and $\varphi_p$, $a_p$ leads to
$$-ihe^{-i\varphi_s/h}B_{\lambda,\mu}\left(e^{i\varphi_s/h}a_s\right)|_{x_1=0}
-ihe^{-i\varphi_p/h}B_{\lambda,\mu}\left(e^{i\varphi_p/h}a_p\right)|_{x_1=0}$$ 
$$=\lambda\langle\rho_s\nu-\beta_0,a_s^0\rangle\nu+\mu\langle\nu,a_s^0\rangle(\rho_s\nu-\beta_0)+
\mu\langle\nu,\rho_s\nu-\beta_0\rangle a_s^0$$
$$+\lambda\langle\rho_p\nu-\beta_0,a_p^0\rangle\nu+\mu\langle\nu,a_p^0\rangle(\rho_p\nu-\beta_0)
+\mu\langle\nu,\rho_p\nu-\beta_0\rangle a_p^0$$
$$-ih(\lambda+\mu)\langle \nu,a_s^1+a_p^1\rangle\nu-ih\mu(a_s^1+a_p^1)$$
$$=(\lambda+\mu)\langle\nu,\rho_sa_s^0+\rho_pa_p^0\rangle\nu-\lambda\langle\beta_0,f\rangle\nu
-\mu\langle\nu,f\rangle\beta_0+\mu(\rho_s a_s^0+\rho_p a_p^0)$$
$$-ih(\lambda+\mu)\langle \nu,a_s^1+a_p^1\rangle\nu-ih\mu(a_s^1+a_p^1)$$
$$=(m_d+hq)\chi f,$$
where
$$q=-i\left(c_s\Pi_s(\nu)+c_p\Pi_p(\nu)\right)\left(A_s^1+A_p^1\right)$$
$$=-i\left(c_s\Pi_s(\nu)+c_p\Pi_p(\nu)\right)U^t(\gamma\widetilde\nabla_{x'}(\rho_p-\rho_s))\Pi_p(\nu)T,$$
   where $\widetilde\nabla_{x'}=(0,\nabla_{x'})$, and 
 $$m_d=\left(c_s\Pi_s(\nu)+c_p\Pi_p(\nu)\right)\left(\rho_sA_s^0+\rho_pA_p^0\right)-\lambda\beta_0\otimes\nu-\mu\nu\otimes\beta_0$$
$$=-\lambda\beta_0\otimes\nu-\mu\nu\otimes\beta_0$$ 
$$ +\left(c_s\Pi_s(\nu)+c_p\Pi_p(\nu)\right)
\left(\rho_sU^t(\rho_s\nu-\beta_0)\Pi_s(\nu)+\rho_pU^t(\rho_p\nu-\beta_0)\Pi_p(\nu)\right)T$$
$$=-\lambda\beta_0\otimes\nu-\mu\nu\otimes\beta_0 
+\left(c_s\rho_s^2\Pi_s(\nu)+c_p\rho_p^2\Pi_p(\nu)\right)T$$
$$-\left(c_s\Pi_s(\nu)+c_p\Pi_p(\nu)\right)U^t(\beta_0)\left(\rho_s\Pi_s(\nu)+\rho_p\Pi_p(\nu)\right)T.$$
 Hence
 \begin{equation}\label{eq:5.24}
 m_d^0(\zeta):=\Lambda m_d\Lambda^{-1}=\lambda\zeta\otimes e_1+\mu e_1\otimes\zeta 
+\left(c_s\rho_s^2\Pi_s(e_1)+c_p\rho_p^2\Pi_p(e_1)\right)T_d$$
$$+\left(c_s\Pi_s(e_1)+c_p\Pi_p(e_1)\right)U^t_0(\zeta)\left(\rho_s\Pi_s(e_1)+\rho_p\Pi_p(e_1)\right)T_d,
\end{equation}
where $T_d(\zeta):=\Lambda T\Lambda^{-1}=W_0^{-1}$.
We will first compute $m_d^0$ when $d=2$. We have
$$U_0^t(\zeta)=\left(
\begin{array}{cc}
0&-\zeta_2\\
\zeta_2&0
\end{array}
\right),$$ 

$$
W_0=\left(
\begin{array}{cc}
\rho_p&-\zeta_2\\
\zeta_2&\rho_s
\end{array}
\right)$$
and hence
$$T_2(\zeta_2)=(r_0+\rho_s\rho_p)^{-1}\left(
\begin{array}{cc}
\rho_s&\zeta_2\\
-\zeta_2&\rho_p
\end{array}
\right).$$
Then we have
$$\lambda\zeta\otimes e_1+\mu e_1\otimes\zeta =\left(
\begin{array}{cc}
0&\lambda\zeta_2\\
\mu\zeta_2&0
\end{array}
\right),$$
 $$\left(c_s\rho_s^2\Pi_s(e_1)+c_p\rho_p^2\Pi_p(e_1)\right)T_2$$
 $$=(r_0+\rho_s\rho_p)^{-1}
 \left(
 \begin{array}{cc}
c_p\rho_p^2&0\\
0&c_s\rho_s^2
\end{array}
\right)
 \left(
\begin{array}{cc}
\rho_s&\zeta_2\\
-\zeta_2&\rho_p
\end{array}
\right)$$ 
$$=(r_0+\rho_s\rho_p)^{-1}
 \left(
 \begin{array}{cc}
c_p\rho_p^2\rho_s&\zeta_2c_p\rho_p^2\\
-\zeta_2c_s\rho_s^2&c_s\rho_s^2\rho_p
\end{array}
\right),$$
$$\left(c_s\Pi_s(e_1)+c_p\Pi_p(e_1)\right)U^t_0(\zeta)\left(\rho_s\Pi_s(e_1)+\rho_p\Pi_p(e_1)\right)T_2$$

$$=(r_0+\rho_s\rho_p)^{-1}
 \left(
 \begin{array}{cc}
0&-\zeta_2c_p\rho_s\\
\zeta_2c_s\rho_p&0
\end{array}
\right)
 \left(
\begin{array}{cc}
\rho_s&\zeta_2\\
-\zeta_2&\rho_p
\end{array}
\right)$$ 
$$=(r_0+\rho_s\rho_p)^{-1}
 \left(
 \begin{array}{cc}
\zeta_2^2c_p\rho_s&-\zeta_2c_p\rho_s\rho_p\\
\zeta_2c_s\rho_s\rho_p&\zeta_2^2c_s\rho_p
\end{array}
\right).$$
Since in this case $\zeta_2=\sqrt{r_0}$, 
an easy computation leads to the formula
 \begin{equation}\label{eq:5.25}
 m_2^0:=(r_0+\rho_s\rho_p)^{-1}\left(
\begin{array}{cc}
z^2n\rho_s&-2\mu\sqrt{r_0}(r_0+\rho_s\rho_p)+z^2n\sqrt{r_0}\\
2\mu\sqrt{r_0}(r_0+\rho_s\rho_p)-z^2n\sqrt{r_0}&z^2n\rho_p
\end{array}
\right).
\end{equation}
Let now $d\ge 3$. Then, in view of (\ref{eq:2.4}) and (\ref{eq:2.6}), we have 
$$\Theta_k(\zeta)^{-1}T_d(\zeta)\Theta_k(\zeta)=\widetilde{T_2}(\sqrt{r_0})+\rho_s^{-1}\sum_{j=3}^de_j\otimes e_j$$
for $\zeta/|\zeta|\in {\cal U}_k$. 
Using this and (\ref{eq:5.24}) one can easily obtain the formula
\begin{equation}\label{eq:5.26}
\Theta_k(\zeta)^{-1}m_d^0(\zeta)\Theta_k(\zeta)=\widetilde{m_2^0}+c_s\rho_s\sum_{j=3}^de_j\otimes e_j=M_d
\end{equation}
for $\zeta/|\zeta|\in {\cal U}_k$. Let $\phi_k\in C^\infty(\mathbb{S}^{d-2})$, $0\le\phi_k\le 1$, $k=1,...,K$, $\sum_{k=1}^K\phi_k=1$,
be a partition of the unity such that supp$\,\phi_k\subset{\cal U}_k$. 
Then we conclude from (\ref{eq:5.25}) and (\ref{eq:5.26}) that 
$$\chi m_d=\sum_{k=1}^K\phi_k(\zeta/|\zeta|)\chi{\cal J}_kM_d{\cal J}_k^{-1}$$
 where 
$${\cal J}_k(x',\xi')=\Lambda(x')^{-1}\Theta_k(\zeta(x',\xi')).$$
\eproof
 
  In what follows we will bound the norm of the difference between the DN map and the operator ${\rm Op}_h(m_d)$. To this end, 
observe that $\widetilde u$ satisfies the equation 
$$(h^2\Delta_{\lambda,\mu}+z^2n)\widetilde u=h\widetilde v,$$
where the function $\widetilde v$ is of the form
$$\widetilde v=(2\pi h)^{-d+1}\int\int e^{\frac{i}{h}\langle y'-x',\xi'\rangle}
\left(e^{i\widetilde\varphi_s/h}B_s+e^{i\widetilde\varphi_p/h}B_p\right)f(y')d\xi'dy'$$ 
$$={\rm Op}_h\left(e^{i\widetilde\varphi_s/h}B_s+e^{i\widetilde\varphi_p/h}B_p\right)f$$
with some matrix-valued functions $B_s$ and $B_p$. To find them we will use the identity
$$e^{-i\varphi/h}(h^2\Delta_{\lambda,\mu}+z^2n(x))\left(e^{i\varphi/h}a\right)
=\left(-P(x,\gamma\nabla_x\varphi)+z^2n\right)a+h^2\Delta_{\lambda,\mu}a+hL(\varphi,A)f,$$
where $a$ is a vector-valued function of the form $a=A(x,\xi')\chi(x',\xi')f(y')$ and $L$ is a matrix-valued function of the form
$$L(\varphi,A)=\sum_{|\alpha|+|\beta|\le 2,\,|\alpha|\ge 1}L_{\alpha,\beta}(x)\partial_x^\alpha\varphi\partial_x^\beta(\chi A),$$
$L_{\alpha,\beta}$ being smooth matrix-valued functions depending only on the variable $x$. Observe also that 
$\Delta_{\lambda,\mu}a=G(A)f$, where $G(A)$ is a matrix-valued function of the form
$$G(A)=\sum_{1\le |\alpha|\le 2}G_{\alpha}(x)\partial_x^\alpha(\chi A).$$
We would like to apply the above identity to $\varphi_s$, $a_s=A_s\chi f$ and $\varphi_p$, $a_p=A_p\chi f$.  
In view of (\ref{eq:5.1}) and (\ref{eq:5.8}), we have 
$$\left(P(x,\gamma\nabla_x\varphi_s)-z^2n\right)a_s=\left(c_s-z^2n(\gamma\nabla_x\varphi_s)^{-2}\right)\Pi_s(\gamma\nabla_x\varphi_s)a_s$$
 $$=x_1^N\Phi_s(\gamma\nabla_x\varphi_s)^{-2}\Pi_s(\gamma\nabla_x\varphi_s)a_s.$$
 By the above identities we get
 $$B_s=h\left[\Delta_{\lambda,\mu},\Psi\right]\chi A_s-h^{-1}x_1^N\Psi\Phi_s(\gamma\nabla_x
 \varphi_s)^{-2}\Pi_s(\gamma\nabla_x\varphi_s)\chi A_s$$
$$+\Psi L(\varphi_s,A_s)+h\Psi G(A_s)$$
and similarly for $B_p$. 
Let $u$ satisfy equation (\ref{eq:1.1}) with $u|_\Gamma={\rm Op}_h(\chi)f$. Then $u-\widetilde u$ satisfies equation (\ref{eq:4.1}) with
$v$ replaced by $\widetilde v$. Therefore, by 
(\ref{eq:4.2}) we get the estimate
\begin{equation}\label{eq:5.27}
\left\|{\cal N}(z,h){\rm Op}_h(\chi)f+ihB_{\lambda,\mu}\widetilde u|_{x_1=0}\right\|_{L^2(\Gamma;\mathbb{C}^d)}\lesssim h^{1/2}\theta^{-1/2}\|\widetilde v\|_{L^2(\Omega;\mathbb{C}^d)}.
\end{equation}
Theorem 1.1 follows from (\ref{eq:5.27}) together with Lemma 5.5 and the following

\begin{lemma} For $N$ big enough depending on $\epsilon$ and $\varepsilon$ we have the estimates
\begin{equation}\label{eq:5.28}
\left\|{\rm Op}_h(\chi q)f\right\|_{L^2(\Gamma;\mathbb{C}^d)}\lesssim \theta^{-1/2-\ell}\|f\|_{H_h^2(\Gamma;\mathbb{C}^d)},
\end{equation}
\begin{equation}\label{eq:5.29}
\|\widetilde v\|_{L^2(\Omega;\mathbb{C}^d)}\lesssim h^{1/2}\theta^{-3/2-\ell}\|f\|_{H_h^3(\Gamma;\mathbb{C}^d)}.
\end{equation}
\end{lemma}

Indeed, we have 
\begin{equation}\label{eq:5.30}
\left\|{\cal N}(z,h){\rm Op}_h(\chi)f-{\rm Op}_h(\chi m_d)f\right\|_{L^2(\Gamma;\mathbb{C}^d)}\lesssim h\theta^{-2-\ell}\|f\|_{H_h^3(\Gamma;\mathbb{C}^d)}.
\end{equation}
We can now take a partition of the unity $\chi_j$, $j=1,...,J$, $0\le\chi_j\le 1$, $\sum_{j=1}^J\chi_j=1$, such that (\ref{eq:5.30}) holds with 
$\chi$ replaced by each $\chi_j$. Moreover, to each $\chi_j$ we can associate a smooth matrix-valued function $\Lambda_j(x')$
such that $\Lambda_j(x')\nu(x')=e_1$ and $\Lambda_j^{-1}(x')=\Lambda_j^t(x')$ in $\pi_{x'}({\rm supp}\,\chi_j)$. 
Thus, summing up all estimates (\ref{eq:5.30}) leads to (\ref{eq:1.4}) with
$$m_d=\sum_{j=1}^J\chi_j m_d=\sum_{j=1}^J\sum_{k=1}^K\chi_j\phi_{k,j}{\cal J}_{k,j}M_d{\cal J}_{k,j}^{-1}$$
where 
$$\phi_{k,j}(x',\xi')=\phi_k\left(-\Lambda_j(x')\beta_0(x',\xi')/\sqrt{r_0(x',\xi')}\right),$$
$${\cal J}_{k,j}(x',\xi')=\Lambda_j(x')^{-1}\Theta_k(-\Lambda_j(x')\beta_0(x',\xi')).$$

\section{Proof of Lemma 5.6}

In view of (\ref{eq:3.2}), we have
$$\nabla_{x'}\rho_s\in S_{2,2}^{-1}(|\rho_s|)+S_{0,1}^1\subset S_{1,1}^{-1/2}(\theta)+S_{0,1}^1$$
and similarly for $\rho_p$. Therefore, we have 
$$U^t(\gamma\widetilde\nabla_{x'}(\rho_p-\rho_s))\in S_{1,1}^{-1/2}(\theta)+S_{0,1}^1$$
which together with (\ref{eq:5.11}) yield
\begin{equation}\label{eq:6.1}
\chi q\in S_{1,1}^{-1/2-\ell}(\theta)+S_{0,1}^2.
\end{equation}
Now (\ref{eq:5.28}) follows from (\ref{eq:6.1}) and Proposition 3.1. Furthermore, it is easy to see that (\ref{eq:5.29}) is
a consequence of the following

\begin{lemma} For $N$ big enough depending on $\epsilon$ and $\varepsilon$ we have the estimate 
\begin{equation}\label{eq:6.2}
\left\|{\rm Op}_h\left(e^{i\widetilde\varphi_s/h}B_s\right)
\right\|_{H_h^3(\Gamma;\mathbb{C}^d)\to L^2(\Gamma;\mathbb{C}^d)}\lesssim h+\theta^{-1-\ell}e^{-Cx_1\theta/h}
\end{equation}
and similarly for $e^{i\widetilde\varphi_p/h}B_p$, where $C>0$ is the same constant as in Lemma 5.2. 
\end{lemma}

Indeed, we have
$$\|\widetilde v\|^2_{L^2(\Omega;\mathbb{C}^d)}\lesssim h^2\|f\|^2_{H_h^3(\Gamma;\mathbb{C}^d)}+\theta^{-2-2\ell}\|f\|^2_{H_h^3(\Gamma;\mathbb{C}^d)}\int_0^\infty
e^{-2Cx_1\theta/h}dx_1$$
$$\lesssim h\theta^{-3-2\ell}\|f\|^2_{H_h^3(\Gamma;\mathbb{C}^d)}.$$
 \eproof
 
 {\it Proof of Lemma 6.1}. Observe first that by (\ref{eq:3.2}) we have
 \begin{equation}\label{eq:6.3} 
\phi_0\left(x_1/|\rho_s|^3\delta\right)\phi_0(x_1\langle\xi'\rangle^\varepsilon/\delta)\in S_{2,2}^0(|\rho_s|)+S_{0,1}^0\subset S_{1,1}^0(\theta)+S_{0,1}^0
\end{equation}
 uniformly in $x_1$, and similarly with $|\rho_s|$ replaced by $|\rho_p|$. In view of the choice of the function $\chi$,
 it is easy to see that (\ref{eq:6.3}) implies
 \begin{equation}\label{eq:6.4} 
\chi\Psi\in S_{1,1}^0(\theta)+S_{0,1}^0
\end{equation}
 uniformly in $x_1$. On the other hand, by Lemma 5.2,
 \begin{equation}\label{eq:6.5} 
e^{Cx_1\theta/h}e^{i\widetilde\varphi_s/h}\in S_{1,1}^0(\theta)+S_{0,1}^0
\end{equation}
on supp$\,\Psi$, uniformly in $x_1$, and similarly with $\widetilde\varphi_s$ replaced by $\widetilde\varphi_p$. By (\ref{eq:6.4}) and (\ref{eq:6.5}), 
\begin{equation}\label{eq:6.6} 
e^{Cx_1\theta/h}e^{i\widetilde\varphi_s/h}\chi\Psi\in S_{1,1}^0(\theta)+S_{0,1}^0
\end{equation}
 uniformly in $x_1$. Furthermore, it is easy to see that Lemma 5.1 yields
 \begin{equation}\label{eq:6.7} 
\partial_{x'}^\alpha\varphi_s\in S_{2,2}^0(|\rho_s|)+S_{0,1}^1\subset S_{1,1}^0(\theta)+S_{0,1}^1,\quad 1\le|\alpha|\le 2,
\end{equation}
\begin{equation}\label{eq:6.8} 
\partial_{x_1}^k\varphi_s\in S_{2,2}^{4-3k}(|\rho_s|)+S_{0,1}^1\subset 
 \left\{
\begin{array}{l}
S_{1,1}^0(\theta)+S_{0,1}^1\quad\mbox{if}\quad k=1,\\
S_{1,1}^{2-3k/2}(\theta)+S_{0,1}^1\quad\mbox{if}\quad k\ge 2,
\end{array}
\right.
\end{equation}
on supp$\,\Psi$, uniformly in $x_1$, and similarly with $\varphi_s$ replaced by $\varphi_p$. By (\ref{eq:6.8}),
\begin{equation}\label{eq:6.9} 
\partial_{x_1}^kU^t(\gamma\nabla_x\varphi_s), \partial_{x_1}^kU^t(\gamma\nabla_x\varphi_p)\in  
 \left\{
\begin{array}{l}
S_{1,1}^0(\theta)+S_{0,1}^1\quad\mbox{if}\quad k=0,\\
S_{1,1}^{(1-3k)/2}(\theta)+S_{0,1}^1\quad\mbox{if}\quad k\ge 1,
\end{array}
\right.
\end{equation}
on supp$\,\Psi$, uniformly in $x_1$. 
By (\ref{eq:6.9}) and (\ref{eq:5.11}),
\begin{equation}\label{eq:6.10} 
\partial_{x_1}^kA_s, \partial_{x_1}^kA_p\in  
 \left\{
\begin{array}{l}
S_{1,1}^{-\ell}(\theta)+S_{0,1}^2\quad\mbox{if}\quad k=0,\\
S_{1,1}^{-\ell+(1-3k)/2}(\theta)+S_{0,1}^2\quad\mbox{if}\quad k\ge 1,
\end{array}
\right.
\end{equation}
on supp$\,\Psi$, uniformly in $x_1$. It is easy to see that (\ref{eq:6.4}), (\ref{eq:6.8}) and (\ref{eq:6.10}) imply
\begin{equation}\label{eq:6.11} 
\Psi L(\varphi_s,A_s), \Psi L(\varphi_p,A_p)\in S_{1,1}^{-1-\ell}(\theta)+S_{0,1}^3,
 \end{equation}
 \begin{equation}\label{eq:6.12} 
\Psi G(A_s), \Psi G(A_p)\in S_{1,1}^{-5/2-\ell}(\theta)+S_{0,1}^2.
 \end{equation}
 By (\ref{eq:6.6}), (\ref{eq:6.11}) and (\ref{eq:6.12}) we conclude
 \begin{equation}\label{eq:6.13} 
e^{Cx_1\theta/h}e^{i\widetilde\varphi_s/h}\Psi\left(L(\varphi_s,A_s)+hG(A_s)\right)\in S_{1,1}^{-1-\ell}(\theta)+S_{0,1}^3
\end{equation}
 as long as $\theta\ge h^{2/5-\epsilon}$. Thus, by (\ref{eq:6.13}) and Proposition 3.1 we obtain
 \begin{equation}\label{eq:6.14}
\left\|{\rm Op}_h\left(e^{i\widetilde\varphi_s/h}\Psi\left(L(\varphi_s,A_s)+hG(A_s)\right)\right)
\right\|_{H_h^3(\Gamma;\mathbb{C}^d)\to L^2(\Gamma;\mathbb{C}^d)}\lesssim \theta^{-1-\ell}e^{-Cx_1\theta/h}.
\end{equation}
 Furthermore, since 
 $$x_1^Ne^{-Cx_1\theta/h}\lesssim h^N\theta^{-N},\quad x_1^Ne^{-Cx_1|\xi'|/h}\lesssim h^N|\xi'|^{-N},$$
 we deduce from Lemma 5.2 that
 \begin{equation}\label{eq:6.15}
 h^{-N}x_1^Ne^{i\widetilde\varphi_s/h}\in S_{1,1}^{-N}(\theta)+S_{0,1}^{-N}
 \end{equation}
 uniformly in $x_1$ and $h$. By (\ref{eq:6.15}) and (\ref{eq:5.3}),
 \begin{equation}\label{eq:6.16}
 h^{-N}x_1^Ne^{i\widetilde\varphi_s/h}\Phi_s\in S_{1,1}^{1-5N/2}(\theta)+S_{0,1}^{2-N}
 \end{equation}
 on supp$\,\Psi$, uniformly in $x_1$ and $h$. On the other hand, it follows from (\ref{eq:6.7}) and (\ref{eq:6.8}) that
 \begin{equation}\label{eq:6.17}
 \Pi_s(\gamma\nabla_x\varphi_s)\in S_{1,1}^0(\theta)+S_{0,1}^2
 \end{equation}
 on supp$\,\Psi$. Taking $N$ big enough, depending on $\varepsilon$, and $\delta$ small enough, we can arrange
 $$\left|z^2n+x_1^N\Phi_s\right|\ge C-x_1^N|\Phi_s|\ge C-{\cal O}(\delta)\ge C/2$$
 on supp$\,\Psi$, with some constant $C>0$. Therefore, using the eikonal equation (\ref{eq:5.1}) we can write
 $$(\gamma\nabla_x\varphi_s)^{-2}=c_s\left(z^2n+x_1^N\Phi_s\right)^{-1}.$$
 In view of (\ref{eq:5.3}) we have
 $$z^2n+x_1^N\Phi_s\in S_{1,1}^0(\theta)+S_{0,1}^0$$
 on supp$\,\Psi$. Thus we obtain
 \begin{equation}\label{eq:6.18}
 (\gamma\nabla_x\varphi_s)^{-2}\in S_{1,1}^0(\theta)+S_{0,1}^0
 \end{equation}
 on supp$\,\Psi$. By (\ref{eq:6.10}), (\ref{eq:6.16}), (\ref{eq:6.17}) and (\ref{eq:6.18}) we conclude
 \begin{equation}\label{eq:6.19}
 h^{-N}x_1^Ne^{i\widetilde\varphi_s/h}\Phi_s\Psi(\gamma\nabla_x\varphi_s)^{-2}\Pi_s(\gamma\nabla_x\varphi_s)\chi A_s\in S_{1,1}^{-5N/2}(\theta)+S_{0,1}^0
 \end{equation}
 provided $N\ge 6$. It follows from (\ref{eq:6.19}) and Proposition 3.1 that
 \begin{equation}\label{eq:6.20}
\left\|{\rm Op}_h\left(h^{-1}x_1^Ne^{i\widetilde\varphi_s/h}\Phi_s\Psi(\gamma\nabla_x\varphi_s)^{-2}\Pi_s(\gamma\nabla_x\varphi_s)\chi A_s\right)
\right\|_{L^2(\Gamma;\mathbb{C}^d)\to L^2(\Gamma;\mathbb{C}^d)}$$
$$\lesssim h^{N-1}\theta^{-5N/2}\lesssim h^{5\epsilon N/2-1}\lesssim h
\end{equation}
  as long as $\theta\ge h^{2/5-\epsilon}$ and $N\ge 4/5\epsilon$. Let $\chi$ be of compact support and suppose that
  supp$\,\chi\cap\Sigma_s\neq\emptyset$, supp$\,\chi\cap\Sigma_p=\emptyset$. Then 
$\left[\Delta_{\lambda,\mu},\Psi\right]\chi=0$ for $x_1\le \delta_1|\rho_s|^3$ for some constant $\delta_1>0$. Therefore, on 
supp$\,\left[\Delta_{\lambda,\mu},\Psi\right]\chi$ we have the bounds
 $$e^{-Cx_1\theta/h}\le e^{-C\delta_1|\rho_s|^3\theta/h}\le e^{-\widetilde C\theta^{5/2}/h}\lesssim h^N\theta^{-5N/2}.$$
 Clearly, we have similar bounds when supp$\,\chi\cap\Sigma_p\neq\emptyset$, supp$\,\chi\cap\Sigma_s=\emptyset$. When  
  supp$\,\chi\cap\Sigma_s=\emptyset$, supp$\,\chi\cap\Sigma_p=\emptyset$, then 
 $\left[\Delta_{\lambda,\mu},\Psi\right]\chi=0$ for $x_1\le \delta_2$ for some constant $\delta_2>0$. So, in this case
 the above bounds still hold. Let now $\chi\in S_{0,1}^0$ be such that ${\rm supp}\,\chi\subset{\rm supp}(1-\eta)$. 
 Then 
 $\left[\Delta_{\lambda,\mu},\Psi\right]\chi=0$ for $x_1\le \delta_3\langle\xi'\rangle^{-\varepsilon}$ for some constant $\delta_3>0$.
 Hence, on 
supp$\,\left[\Delta_{\lambda,\mu},\Psi\right]\chi$ we have the bounds
  $$e^{-Cx_1|\xi'|/h}\le e^{-\widetilde C|\xi'|^{1-\varepsilon}/h}\lesssim h^N|\xi'|^{-N(1-\varepsilon)}.$$
  Therefore, by Lemma 5.2 and (\ref{eq:6.10}) we get
  \begin{equation}\label{eq:6.21}
 h^{-N}e^{i\widetilde\varphi_s/h}\left[\Delta_{\lambda,\mu},\Psi\right]\chi A_s\in S_{1,1}^{-\ell-1-5N/2}(\theta)+S_{0,1}^0
 \end{equation}
 for $N$ big enough. 
 It follows from (\ref{eq:6.21}) and Proposition 3.1 that
 \begin{equation}\label{eq:6.22}
\left\|{\rm Op}_h\left(he^{i\widetilde\varphi_s/h}\left[\Delta_{\lambda,\mu},\Psi\right]\chi A_s\right)
\right\|_{L^2(\Gamma;\mathbb{C}^d)\to L^2(\Gamma;\mathbb{C}^d)}$$
$$\lesssim h^{N+1}\theta^{-\ell-1-5N/2}\lesssim h^{5\epsilon N/2+1-(\ell+1)(2/5-\epsilon)}\lesssim h
\end{equation}
  as long as $\theta\ge h^{2/5-\epsilon}$ and $N$ big enough. Now the estimate (\ref{eq:6.2}) follows from
  (\ref{eq:6.14}), (\ref{eq:6.20}) and (\ref{eq:6.22}).
\eproof

\end{document}